\documentclass{amsart}
\usepackage[T1]{fontenc}   
\usepackage{graphicx}
\usepackage{mathptmx}
\usepackage{hyperref}
\usepackage[all,cmtip]{xy}

\usepackage{amsmath,amssymb,amscd,amsthm,latexsym}
\usepackage[applemac]{inputenc}

\usepackage{latexsym,enumerate,hyperref,amsmath,amssymb,
amscd,enumerate,engord,amsfonts,graphics}

\newtheorem{theorem}{Theorem}[section]
\newtheorem{lemma}[theorem]{Lemma}
\newtheorem{proposition}[theorem]{Proposition}
\theoremstyle{definition}
\newtheorem{definition}[theorem]{Definition}

\newtheorem{remark}[theorem]{Remark}
\newtheorem{example}[theorem]{Example}
\numberwithin{equation}{section}

\def\R{{\mathbb{R}}}
\def\Z{{\mathbb{Z}}}
\def\P{{\mathbb{P}}}
\def\XX{{\mathcal{P}}}

\DeclareMathOperator{\TC}{{\sf TC}}
\DeclareMathOperator{\cat}{\mathrm{cat}}

\DeclareMathOperator{\ev}{\mathrm{ev}}
\DeclareMathOperator{\id}{\mathrm{id}}
\DeclareMathOperator{\zcl}{\mathrm{zcl}}

\begin{document}

\renewcommand{\bf}{\bfseries}
\renewcommand{\sc}{\scshape}

\title[Motion planning in connected sums of projective spaces]%
{Motion planning in connected sums of real projective spaces}

\author{Daniel C. Cohen}
\address{Department of Mathematics, Louisiana State University, Baton Rouge, LA 70803 USA;
\quad  {\rm{\texttt{\href{mailto:cohen@math.lsu.edu}{cohen@math.lsu.edu}}; \quad \texttt{\href{http://www.math.lsu.edu/~cohen/}{www.math.lsu.edu/\char'176cohen}}}}}
\thanks{The first author was partially supported by the Simons Foundation and by the Mathematisches Forschungsinstitut Oberwolfach.}

\author{Lucile Vandembroucq}
\address{Centro de Matem\'atica, Universidade do Minho, Campus de Gualtar, 4710-057 Braga, Portugal; \quad
{\rm{\texttt{\href{mailto:lucile@math.uminho.pt}{lucile@math.uminho.pt}}}}}
\thanks{The second author was partially supported by FCT-UID/MAT/00013/2013 and by the Polish National Science Centre grant $2016/21/P/ST1/03460$ within the European Union's Horizon 2020 research and innovation programme under the Marie Sk\l odowska-Curie grant agreement No. 665778.  }



\begin{abstract} 
The topological complexity $\TC(X)$ is a homotopy invariant of a topological space $X$, motivated by robotics, and providing a measure of the navigational complexity of $X$. The topological complexity of a connected sum of real projective planes, that is, a high genus nonorientable surface, is known to be maximal. We use algebraic tools to show that the analogous result holds for connected sums of higher dimensional real projective spaces.
\end{abstract}

\maketitle

\section{\bf Introduction}
Let $X$ be a finite, path-connected CW-complex.  Viewing $X$ as the space of configurations of a mechanical system, the motion planning problem consists of constructing an algorithm which takes as input pairs of configurations $(x_0,x_1) \in X \times X$, and produces a continuous path $\gamma\colon [0,1] \to X$ from the initial configuration $x_0=\gamma(0)$ to the terminal configuration $x_1=\gamma(1)$.  The motion planning problem is of significant interest in robotics, see, for example, Latombe~\cite{La} and Sharir \cite{Sh}.

A topological approach to the motion planning problem is developed by Farber, see  \cite{Far,FarRobots,FarInv}. Let $I=[0,1]$ be the unit interval, and let $X^I$ be the space of continuous paths $\gamma\colon I \to X$ (with the compact-open topology). The map $\ev \colon X^I \to X \times X$ defined by sending a path to its endpoints, $\ev(\gamma)=(\gamma(0),\gamma(1))$, is a fibration, with fiber $\Omega(X)$, the based loop space of $X$. The motion planning problem requests a section of this fibration, a map $s \colon X \times X \to X^I$ satisfying 
$\ev \circ s = \id_{X\times X}$. It would be desirable for the motion planning algorithm to depend continuously on the input.  However, there exists a globally continuous section $s\colon X\times X \to PX$ if and only if $X$ is contractible, see \cite[Thm.~1]{Far}. This prompts the study of the discontinuities of such algorithms, and leads to the following definition from \cite{FarRobots}.

\begin{definition} \label{def:motion planner}
A \emph{motion planner} for $X$ is a collection of subsets $F_0,F_1,\dots,F_m$ of $X\times X$ and continuous maps $s_i\colon F_i \to PX$ such that
\begin{enumerate}
\item the sets $F_i$ are pairwise disjoint, $F_i \cap F_j=\emptyset$ if $i \neq j$, and cover $X \times X$, 
\[
X\times X = F_0 \cup F_1 \cup \cdots \cup F_m;
\]
\item $\ev \circ s_i = \id_{F_i}$ for each $i$; and
\item each $F_i$ is a Euclidean neighborhood retract.
\end{enumerate}
\end{definition}
\noindent Refer to the sets $F_i$ as local domains of the motion planner, and the maps $s_i$ as local rules.  Call a motion planner optimal if it requires a minimal number of local domains (resp., rules).  

\begin{definition} \label{def:TC} For a finite, path-connected CW-complex $X$, the (reduced) \emph{topological complexity} of $X$, $\TC(X)$, is one less than the number of local domains in an optimal motion planner for $X$, $\TC(X)=m$ if there exists an optimal motion planner $F_0,F_1,\dots,F_m$ for $X$.  
\end{definition}

\subsection{\bf Motion planning in cell complexes}\label{sec:mp}
We briefly recall from \cite[\S3]{FarRobots} a construction of a motion planner for a finite cell complex. Recall that $X$ is a finite, path-connected CW-complex, and let $X^k$ be the $k$-dimensional skeleton of $X$.  Assume that $\dim(X)=n$, and for $k=0,1,\dots,n$, let $V^k = X^k \setminus X^{k-1}$ be the union of the open $k$-cells of $X$. For $i=0,1,\dots,2n$, the sets 
$
F_i = \bigcup_{k+l=i} V^k \times V^l \subset X\times X
$
are homeomorphic to disjoint unions of balls, so are Euclidean neighborhood retracts. Note that $F_0 \cup F_1 \cup \cdots \cup F_m=X\times X$.

To define a local rule $s_i \colon F_i \to X^I$, since $F_i$ is the union of disjoint sets $V^k \times V^l$ (which are both open and closed in $F_i$), it suffices to construct a continuous map $s_{k,l}\colon V^k \times V^l \to X^I$ satisfying $\ev \circ s_{k,l}=\id_{V^k \times V^l}$. Pick a point $v_k \in V^k$ for each $k$, and fix a path $\gamma_{k,l}$ in $X$ from $v_k$ to $v_l$ for each $k,l$. Then, for any $(x,y) \in V^k\times V^l$, one can construct a path $s_{k,l}(x,y)$ from $x$ to $y$ by first moving from $x$ to $v_k$ in the cell $V^k$, then traversing the fixed path  $\gamma_{k,l}$, and finally moving from $v_l$ to $y$ in $V^l$.

This construction exhibits a motion planner for $X$ with $2\dim(X)+1$ local domains. Consequently, we have the upper bound $\TC(X) \le 2 \dim(X)$ (for a finite, path connected CW-complex $X$). This upper bound is achieved by many spaces of interest in topology and applications. For instance, it is well known that $\TC(\Sigma_g)=4$ for an orientable surface $\Sigma_g$ of genus $g\ge 2$, see \cite{Far}. More recent work of Dranishnikov \cite{Dranishnikov,Dranishnikov2} and the authors \cite{CV} shows that the same holds for nonorientable surfaces of high genus. Observe that the construction above provides an optimal motion planner in these instances.

\subsection{\bf Main result}
The objective of this note is to establish a higher dimensional analog of these last results. 
Let $\XX^n_g=\R\P^n\#\cdots \# \R\P^n$ be the connected sum of $g$ copies of the real projective space $\R\P^n$. 
\begin{theorem} \label{main} For $n\geq 2$ and $g\geq 2$, we have $\TC(\XX_g^n)=2n$.
\end{theorem}

Thus, applying the construction in \S\ref{sec:mp} above to a standard CW decomposition of the space $\XX_g^n$ yields an optimal motion planner for this space.

When $n=2$, $\XX_g^2=N_g$ is the nonorientable surface of genus $g$, and it has been established in \cite{CV} that $\TC(N_g)=4$ for $g\geq 2$, completing results obtained by Dranishnikov \cite{Dranishnikov,Dranishnikov2} in the case $g\geq 4$. 
So we focus on the case $n\ge 3$ below. 
As we will see, the methods developed in \cite{CV} admit extensions to this higher dimensional case.

\begin{remark} 
The case $g=1$, with $\XX_1^n=\R\P^n$, is significantly more subtle. As shown by Farber-Tabachnikov-Yuzvinsky \cite{FTY}, for $n \neq 1,3,7$, the topological complexity and immersion dimension of $\R\P^n$ are equal, $\TC(\R\P^n)=\mathrm{imm}(\R\P^n)$.
\end{remark}

\section{\bf Preliminaries}
Let $p\colon E \to B$ be a fibration. The (reduced) sectional category, or Schwarz genus, of $p$, denoted by $\mathrm{secat}(p)$, is the smallest integer $m$ such that $B$ can be covered by $m+1$ open subsets, over each of which $p$ has a continuous section. Classical references include Schwarz \cite{Sch} and James \cite{James}. The following result makes clear the topological nature of the motion planning problem.

\begin{theorem}[{\cite[cf.\,\S4.2]{FarInv}}] If $X$ is a finite CW-complex, then the topological complexity of $X$ is equal to the sectional category of the path-space fibration $\ev\colon X^I \to X \times X$, $\TC(X)=\mathrm{secat}(\ev)$.
\end{theorem}

The equality $\TC(X)=\mathrm{secat}(\ev\colon X^I \to X \times X)$ yields the following estimates:
\begin{equation*}
\max\{\cat(X), \zcl_{\Bbbk}(X)\} \leq \TC(X) \leq 2\cat(X)\leq 2\dim(X), \ \text{see \cite{Far}}.
\end{equation*}
Here, $\cat(X)$ is the reduced Lusternik-Schnirelmann category of $X$ 
 and $\zcl_{\Bbbk}(X)$ is the {zero-divisors cup-length} of the cohomology of $X$ with coefficients in a field $\Bbbk$. More precisely, $\zcl_{\Bbbk}(X)$ is the nilpotency of the kernel of the cup product $H^*(X;\Bbbk)\otimes H^*(X;\Bbbk)\to H^*(X;\Bbbk)$, the smallest nonnegative integer $n$ such that any $(n+1)$-fold cup product in this kernel is trivial. 
 
As noted in \S\ref{sec:mp}, the upper bound $\TC(X) \le 2 \dim(X)$ may also be obtained from an explicit motion planner construction. 
We will not make further use of the lower bounds $\cat(X)$ and $\zcl_{\Bbbk}(X)$, which are included here primarily for context, and are both insufficient for our purposes. Indeed for $g\geq 2$, one can show that $\cat(\XX^n_g)=n$ and $\zcl_{\Z_2}(\XX^n_g)=2n-1$. 
Following \cite{CV}, we will instead utilize the topological complexity analog of the classical Berstein-Schwarz cohomology class, which informs on the LS category, 
see \cite[Thm. 2.51]{CLOT}.

Let $X$ be a space and $\pi=\pi_1(X)$ its fundamental group.  Let $\Z[\pi]$ be the group ring of $\pi$, $\epsilon\colon \Z[\pi]\to\Z$ the augmentation map, and $I(\pi)=\ker(\varepsilon\colon\Z[\pi]\to \Z)$ the augmentation ideal. Recall that  $\Z[\pi]$ and $I(\pi)$ are both (left) $\Z[\pi\times \pi]$-modules through the action given by:
\[
(a, b) \cdot \sum n_ia_i= \sum n_i (aa_i\bar{b}).
\]
Here $n_i\in \Z$, $a,b,a_i\in \pi$, and $\bar b$ is the inverse of $b$. In general (see \cite[\S6]{Wh}), left $\Z[\pi\times \pi]$-modules correspond to local coefficient systems on $X\times X$, which we denote by the same symbols.

Let ${\mathfrak{v}}={\mathfrak{v}}_X\in H^ 1(X\times X;I(\pi))$ be the Costa-Farber canonical class of $X$ introduced in \cite{CostaFarber}, corresponding to the crossed homomorphism $\pi\times \pi \to I(\pi)$, $(a,b)\mapsto a\bar{b}-1$. The significance of this cohomology class in the context of topological complexity is given by the following result.
\begin{theorem}[{\cite[Theorem 7]{CostaFarber}}]  Suppose that $X$ is a CW-complex of dimension $n\geq 2$. Then $\TC(X)=2n$ if and only if the $2n$-th power of ${\mathfrak{v}}$ does not vanish:
\[
\TC(X)= 2n \Longleftrightarrow  {\mathfrak{v}}^{2n}\neq 0\ \text{in}\   H^{2n}(X\times X;I(\pi)^{\otimes 2n}).
\]
\end{theorem}
\noindent Here $I(\pi)^{\otimes 2n}=I(\pi)\otimes_\Z I(\pi) \otimes_\Z \cdots \otimes_\Z I(\pi)$ is the tensor product of $2n$ copies of $I(\pi)$, with the diagonal action of $\pi\times \pi$.

\section{\bf Reduction to the case $g=2$}

Let $\pi_g$ denote the fundamental group of the space $\XX_g^n$. Since $n\geq 3$, we have $\pi_g=\Z_2\ast \cdots \ast\Z_2$ ($g$ copies).
As in \cite{CV}, we will prove that $\TC(\XX_g^n)=2n$ by proving that the evaluation of $\mathfrak{v}^{2n}\in H^{2n}(\XX_g^n\times \XX_g^n; I(\pi_g)^{\otimes 2n})$ on the $\Z_2$ top class $[\XX_g^n\times \XX_g^n]_{\Z_2}\in H_{2n}(\XX_g^n\times \XX_g^n;\Z_2)$
does not vanish and use the bar resolution to carry out the calculation. As noted in \cite[Corollary 8]{CostaFarber}, if $f:X\to Y=K(\pi,1)$ induces an isomorphism of fundamental groups, we have $(f\times f)^*\mathfrak{v}_Y=\mathfrak{v}_X$. Additionally, as used in \cite{CV} and recalled below, the class $\mathfrak{v}_Y\in H^1(\pi\times \pi;I(\pi))$ and its powers can be represented by the powers of an explicit cocycle $\nu$ defined on the bar resolution of $\pi\times \pi$. Denoting by $f_g\colon \XX_g^n \to K(\pi_g,1)$ the canonical map, that is, the unique (up to homotopy) map such that $\pi_1(f_g)=\id$, we then analyze $\mathfrak{v}^{2n}([\XX_g^n\times \XX_g^n]_{\Z_2})$ using the following composite
\[ H_{2n}(\XX_g^n\times \XX_g^n;\Z_2)\xrightarrow{(f_g\times f_g)_*} H_{2n}(\pi_g\times \pi_g;\Z_2)\xrightarrow{\nu^{2n}_*} I(\pi_g;\Z_2)^ {\otimes 2n}_{\pi_g\times \pi_g}
\]
where $I(\pi_g;\Z_2)=I(\pi_g)\otimes \Z_2$ and $I(\pi_g;\Z_2)^ {\otimes 2n}_{\pi_g\times \pi_g}$ denotes the coinvariants of $I(\pi_g;\Z_2)^ {\otimes 2n}$ with respect to the diagonal action of $\pi_g\times \pi_g$.

As in \cite[Theorem 14]{CV}, the study of the general case $g\geq 2$ can be reduced to the case $g=2$. Consider the projection $\XX_g^n\to \XX_{g-1}^n$ that collapses the last $\R\P^ n$ component of $\XX_g^n$ and induces the projection $\pi_g\to \pi_{g-1}$ which sends the last $\Z_2$ to $1$. We have a (homotopy) commutative diagram
$$
\xymatrixcolsep{5pc}
\xymatrix{
\XX_g^n \ar[d]\ar[r]^{f_g} & K(\pi_g, 1)\ar[d]\\
\XX_{g-1}^n \ar[r]^{f_{g-1}} & K(\pi_{g-1}, 1)
}$$
Since the projection $\XX_g^n\to \XX_{g-1}^n$ fits in a cofibration sequence $\R\P^ {n-1}\to \XX_g^n\to \XX_{g-1}^n$, the induced morphism $H_n(\XX_g^n;\Z_2)\to H_n(\XX_{g-1}^n;\Z_2)$ is an isomorphism. Considering the morphism 
$I(\pi_g;\Z_2)\to I(\pi_{g-1};\Z_2)$ induced by the projection $\pi_g\to \pi_{g-1}$, we then obtain the following commutative diagram
$$\xymatrix{
\Z_2\cong H_{2n}(\XX_g^n\times \XX_g^n;\Z_2) \ar[d]_{\cong}\ar[rrr]^-{\nu_*^{2n}(f_g\times f_g)_*} &&& I(\pi_g;\Z_2)^{\otimes 2n}_{\pi_g\times \pi_g}\ar[d]\\
\Z_2\cong H_{2n}(\XX_{g-1}^n\times \XX_{g-1}^n;\Z_2)  \ar[rrr]^-{\nu_{*}^{2n}(f_{g-1}\times f_{g-1})_*} &&& I(\pi_{g-1};\Z_2)^{\otimes 2n}_{\pi_{g-1}\times \pi_{g-1}}.\\
}$$
in which the left hand vertical map is an isomorphism. Therefore, if the bottom horizontal map does not annihilate the generator, then neither does the top horizontal map. In other words, as in \cite{CV}, the calculation can be reduced to the ``genus'' $g=2$ case. Thus, for $n\geq 3$, Theorem \ref{main}  will follow from
the following proposition which will be proved in the next section.

\begin{proposition}\label{g=2}
For $n\geq 3$, $\mathfrak{v}^{2n}([\XX_2^n\times \XX_2^n]_{\Z_2})\neq 0$.
\end{proposition}

We note that, for $M,N$ closed $n$ dimensional manifolds, a similar argument to the one above permits one to conclude that $\TC(M\#N)=\TC(M)=2n$ as soon as $\mathfrak{v}^{2n}([M\times M]_{\Z_2})$ is non zero. Actually, using \cite[Lemma 7]{Dranishnikov-Sadykov} (and $\Z$-fundamental classes instead of $\Z_2$ top classes), we can see that $\TC(M\#N)$ is maximal as soon as $\TC(M)$ is maximal whenever $N$ is orientable. Note also that, for simply-connected orientable manifolds, Dranishnikov and Sadykov \cite{Dranishnikov-Sadykov} established the more general result that  $\TC(M\#N)\geq\TC(M)$. 

\section{\bf The case $g=2$}

\subsection{\bf Algebraic preliminaries}
Refer to Brown \cite{Brown} and Weibel \cite{Weibel} as standard references for cohomology of groups and homological algebra. We will use the normalized bar resolution $\bar B_{*}(\pi)$ of $\Z$ as a trivial $\Z[\pi]$-module:
\[
\cdots \longrightarrow \bar B_n(\pi)\xrightarrow{\ \partial_n\ } \cdots \longrightarrow \bar B_1(\pi) \xrightarrow{\ \partial_1\ } \bar B_0(\pi)=\Z[\pi] \xrightarrow{\ \varepsilon\ } \Z \longrightarrow 0.
\]
Here $\bar B_n(\pi)$ is the free $\Z[\pi]$-module with basis $\{[g_1|\cdots |g_n], (g_1,\dots ,g_n)\in \bar\pi^n\}$,  where $\bar\pi=\{g\in \pi \mid g\neq 1\}$ and $\partial_n$ is the $\Z[\pi]$ morphism given by 
\begin{multline*}
\partial_n([g_1|\cdots |g_n])=g_1\cdot [g_2|\cdots |g_n]+\sum\limits_{i=1}^{n-1}(-1)^i[g_1|\cdots |g_{i-1}|g_ig_{i+1}|g_{i+2}|\cdots |g_n]\\ + (-1)^n[g_1|\cdots |g_{n-1}]
\end{multline*}
(with $[h_1|\cdots |h_k]=0$ if $h_i=1$ for some $i$). The homology of the space $K(\pi,1)$ (or of the group $\pi$) with coefficients in $\Z_2$ is then the homology of the chain complex $\bar B_{*}(\pi;\Z_2):= \bar B_{*}(\pi)\otimes_{\pi}\Z_2=(\bar B_{*}(\pi))_\pi\otimes\Z_2$ (with differential $\partial \otimes\, \mathrm{id}$). 

We now describe a cycle representing the image of the $\Z_2$ top class of $\XX_2^n=\R\P^n\#\R\P^n$ under the map induced by $f_2:\XX_2^n \to K(\pi_2,1)$. We have $H_i(\pi_2;\Z_2)=H_i(\R\P^{\infty}\vee \R\P^{\infty};\Z_2)$. 
Let $\mathbf{a}_i, \mathbf{b}_i$ be the homology classes (with $\mathbf{a}_0= \mathbf{b}_0$) corresponding to the two branches of the wedge.
As the two projections $\XX_2^n=\R\P^ n\#\R\P^ n \to \R\P^ n$ each induce an isomorphism $H_n(\XX_2^n;\Z_2) \to H_n(\R\P^ n;\Z_2)$, the image of the $\Z_2$ top cell of $\R\P^n\#\R\P^n$ under the map $f_2:\XX_2^n \to K(\pi_2,1)$ can be identified with the element $\mathbf{c}_n=\mathbf{a}_n+ \mathbf{b}_n$ of $H_n(\pi_2;\Z_2)$ and we are reduced to describe cycles representing the classes  $\mathbf{a}_n$ and $\mathbf{b}_n$. 

Writing $\pi_2=\Z_2 * \Z_2=\langle a,b \,|\, a^ 2=1, b^ 2=1\rangle$, the classes  $\mathbf{a}_i, \mathbf{b}_i$ are represented by the following cycles of $B_i(\pi_2;\Z_2)$:
\[ \alpha_i=[a|a|\cdots |a], \qquad \beta_i=[b|b|\cdots |b].
\] 
As our calculation will use portions of the calculation carried out in \cite{CV}, we will use the isomorphism from $\pi_2=\langle a, b ~|~a^2=1,b^2=1\rangle$ to the infinite dihedral group  $D=\langle x, y ~|~yxy=x, x^2=1\rangle$ given by $a\mapsto x$ and $b\mapsto yx$.
We will then work with the following cycles of $B_i(D;\Z_2)$ as representatives of the classes $\mathbf{a}_i,\mathbf{b}_i$:
\[ \alpha'_i=[x|x|\cdots |x],  \qquad \beta'_i=[yx|yx|\cdots |yx].
\]

For $X=K(\pi,1)$, the Costa-Farber $\TC$ canonical class ${\mathfrak v}\in H^ 1(X\times X;I(\pi))$ can be described as the class of the canonical degree $1$ cocycle, 
$\nu\colon \bar B_1(\pi\times \pi)\to I(\pi)$, which is well-defined on the normalized bar resolution and given by
\[
\nu( [(g, h)])= g\bar{h}-1
\]
for $[(g,h)]\in \bar B_1(\pi\times \pi)$, and $\bar{h}=h^{-1}$ as above. As in \cite{CV}, we have the following explicit expression of the $n$-th power of ${\mathfrak v}\in H^ 1(X\times X;I(\pi))$:

\begin{lemma}\label{explicitpower} The $n$-th power of the canonical $\TC$ cohomology class $\mathfrak{v}$ is the class of the cocyle $\nu^n$ of degree $n$ given by
{
\setlength\arraycolsep{2pt}
\begin{eqnarray*}
\nu^n\colon \bar B_n(\pi\times \pi) &\to & I(\pi)^{\otimes n}\\
\left[(g_1,h_1)|\cdots |(g_n, h_n)\right] & \mapsto & \xi\cdot \bigl(u_1-1\bigr)\otimes
\cdots \otimes (g_1\cdots g_{n-1})\bigl(u_n-1\bigr)(\bar{h}_{n-1}\cdots  \bar{h}_{1}),\nonumber
\end{eqnarray*}
}
where $\xi=(-1)^{n(n-1)/2}$ and $u_i=g_i\bar{h}_i$ for each $i$, $1\le i\le n$.
\end{lemma}

We will also use the Eilenberg-Zilber chain equivalence (well-defined on normalized bar resolutions)
\begin{equation} \label{eq:EZ}
EZ\colon \bar B_*(\pi)\otimes \bar B_*(\pi) \longrightarrow \bar B_*(\pi \times \pi),
\end{equation}
which is the $\Z[\pi\times \pi]\cong\Z[\pi]\otimes \Z[\pi]$ morphism given by
\[
\begin{array}{rcl}
EZ_n\colon \bigoplus\limits_{i=0}^{n}\bar B_i(\pi)\otimes \bar B_{n-i}(\pi)& \to & \bar B_n(\pi\times \pi) \\
{[}g_1|\cdots |g_i]\otimes [h_{i+1}|\cdots |h_n] &\mapsto & \sum\limits_{\sigma \in {\mathcal S}_{i,n-i}}\mathrm{sgn}(\sigma)[q_{\sigma^{-1}(1)}| \cdots | q_{\sigma^{-1}(n)}]
\end{array}
\]
where ${\mathcal S}_{i,n-i}$ denotes the set of $(i,n-i)$ shuffles, $\mathrm{sgn}(\sigma)$ is the signature of the shuffle $\sigma$ (which can be omitted over $\Z_2$), and
\[
q_{k}=\begin{cases}
(g_{k}, 1) & \mbox{ if } 1\leq k \leq i, \\
(1, h_{k}) & \mbox{ if } i+1\leq k \leq n.
\end{cases}
\]

\begin{example}\label{example}\rm{We now calculate the explicit expression of $\nu^ 4(EZ(\alpha'_2 \otimes \beta'_2))$ in $I(D;\Z_2)^{\otimes 4}$, which will be useful in the proof of Proposition \ref{g=2}. Since $\alpha'_2=[x|x]$ and $\beta'_2=[yx|yx]$ we have
$$\arraycolsep=1.4pt\def\arraystretch{1.2}
EZ(\alpha'_2\otimes \beta'_2)=\left\{\begin{array}{c}
[x_1|x_1|y_2x_2|y_2x_2]+[x_1|y_2x_2|x_1|y_2x_2]+[x_1|y_2x_2|y_2x_2|x_1]\\
+[y_2x_2|x_1|x_1|y_2x_2]+[y_2x_2|x_1|y_2x_2|x_1]+[y_2x_2|y_2x_2|x_1|x_1]\\
\end{array}\right.,
$$
where $x_1=(x,1)$, $x_2=(1,x)$, $y_1=(y,1)$, $y_2=(1,y)$. 
Using Lemma \ref{explicitpower} together with the fact that $x^ 2=1$ and $(yx)^ 2=1$, we obtain: 
\begin{equation}
\label{expressionexample}
\arraycolsep=1.4pt\def\arraystretch{1.2}
\begin{array}{rcl}
\nu^4(EZ(\alpha'_2\otimes \beta'_2))&=&\left\{\begin{array}{c}
(x-1)\otimes (1-x) \otimes (yx-1) \otimes (1-yx) \\
+(x-1)\otimes x(yx-1) \otimes (1-x)yx \otimes (1-yx) \\
+(x-1)\otimes x(yx-1) \otimes x(1-yx) \otimes (1-{x}) \\
+(y{x}-1)\otimes (x-1)yx \otimes (1-x)yx \otimes (1-yx) \\
+(y{x}-1)\otimes (x-1)yx \otimes x(1-yx) \otimes (1-{x}) \\
+(y{x}-1)\otimes (1-y{x}) \otimes ({x}-1) \otimes (1-{x}) \\
\end{array}\right.
\end{array}
\end{equation}
The image of this expression in the coinvariants $I(D;\Z_2)^{\otimes 4}_{D\times D}$ corresponds to the element $\mathfrak{v}^4(\mathbf{a}_2\times \mathbf{b}_2)\in H_0(D\times D; I(D;\Z_2)^{\otimes 4})\cong I(D;\Z_2)^{\otimes 4}_{D\times D}$.

Let $Y=\langle y \,|\, y^2=1\rangle$ and $Z=\langle z \,|\, z^2=1\rangle$. We have $I(Y;\Z_2)\cong \Z_2(y-1)$ and $I(Z;\Z_2)\cong \Z_2(z-1)$.
Consider the projection $I(D;\Z_2)\to  I(Y;\Z_2)$ sending $x$ to $1$ and the projection $I(D;\Z_2)\to  I(Z;\Z_2)$ sending both $x$ and $y$ to $z$ (and hence $yx\mapsto 1$). After projection onto $I(Y;\Z_2)^{\otimes 2} \otimes I(Z;\Z_2)^{\otimes 2}\cong \Z_2$, one can check that \eqref{expressionexample} yields a unique non-zero term 
$(y-1)\otimes (y-1) \otimes (z-1)\otimes (z-1)$ 
which corresponds to the element $[y_2x_2|y_2x_2|x_1|x_1]\in \bar{B}_4(D\times D)$.
}
\end{example}

\subsection{\bf Proof of Proposition \ref{g=2}}

The statement will follows from the fact that the class $\mathbf{c}_n=(f_2)_*([\XX_2^n]_{\Z_2})\in H_n(\pi_2;\Z_2)=H_n(D;\Z_2)$ corresponding to the $\Z_2$ top class of $\XX_2^n$ satisfies $\mathfrak{v}^{2n}(\mathbf{c}_n \times \mathbf{c}_n)=\mathfrak{v}^{2n}\cap(\mathbf{c}_n \times \mathbf{c}_n) \neq 0$. 

Let $G=D\times D$, $X=K(G,1)$, and let $\Delta\colon X \to X \times X$ be the diagonal map. For the $G$-modules $M=I(D)^{\otimes 2n}$, $M'=I(D)^{\otimes 4}$, and $M''=I(D)^{\otimes 2n-4}$, and the homology and cohomology classes $x\in H_{2n}(G;\Z_2)$ and $y \in H^{2n}(G\times G;M)=H^{2n}(G\times G;M'\otimes M'')$, we have
\[
\Delta_*(x \cap \Delta^*(y)) = \Delta_*(x) \cap y
\]
in $H_0(G \times G; M \otimes \Z_2)$ (cf. \cite[V.10]{Bredon}, see also \cite[Prop. 2.4.2]{Dranishnikov}).
Fixing $y={\mathfrak v}^4 \times {\mathfrak v}^{2n-4} \in H^{2n}(G\times G;M' \otimes M'')$, so that $\Delta^*(y)=
{\mathfrak v}^4 \cup {\mathfrak v}^{2n-4}={\mathfrak v}^{2n}$, we obtain the commuting diagram
\begin{equation*}
\xymatrixcolsep{5pc}
\xymatrix{
H_{2n}(G,\Z_2) \ar[d]_{{\mathfrak v}^{2n}} \ar[r]^-{\Delta_*} & 
H_{2n}(G\times G,\Z_2) \ar[d]^{{\mathfrak v}^4 \times {\mathfrak v}^{2n-4}} \\
H_0(G; M \otimes \Z_2)  \ar[r]^-{\Delta_*} & 
H_0(G\times G; M' \otimes M'' \otimes \Z_2)
}
\end{equation*}
where the vertical maps are cap products with the indicated cohomology classes.

Let $\kappa_{4,2n-4}$ denote the composition of the K\"unneth isomorphism and the projection indicated below
\[
H_{2n}(G\times G;\Z_2) \rightarrow \bigoplus_{i+j=2n} H_i(G;\Z_2) \otimes H_j(G;\Z_2) \twoheadrightarrow H_4(G;\Z_2) \otimes H_{2n-4}(G;\Z_2),
\]
and let $\Delta_{4,2n-4}=\kappa_{4,2n-4} \circ \Delta_*$. Recalling that $G=D\times D$, $M=M'\otimes M''=I(D)^{\otimes 2n}$, and identifying zero-dimensional homology groups, the above commuting diagram yields
\begin{equation}\label{decomposition}
\xymatrixcolsep{5pc}
\xymatrix{
H_{2n}(D\times D,\Z_2) \ar[d]_{{\mathfrak v}^{2n}} \ar[r]^-{\Delta_{4,2n-4}} & 
H_4(D\times D;\Z_2)\otimes H_{2n-4}(D\times D;\Z_2)\ar[d]^{{\mathfrak v}^4\otimes {\mathfrak v}^{2n-4}}\\
I(D;\Z_2)^{\otimes 2n}_{D\times D}
 \ar[r] & I(D;\Z_2)^{\otimes 4}_{D\times D}\otimes I(D;\Z_2)^{\otimes {2n-4}}_{D\times D}
}
\end{equation}

As above, we consider $Y=\langle y \,|\, y^2=1\rangle$ and $Z=\langle z \,|\, z^2=1\rangle$ and the projections $I(D;\Z_2)\to  I(Y;\Z_2)$ and $I(D;\Z_2)\to  I(Z;\Z_2)$.
We then compose the ${\mathfrak v}^{2n-4}$ portion of the right-hand vertical map in the diagram \eqref{decomposition} with the projection 
\[
I(D;\Z_2)^{\otimes 2n-4}_{D\times D} \longrightarrow I(Y;\Z_2)^{\otimes n-2}\otimes I(Z;\Z_2)^{\otimes n-2}.
\]
Observe that 
\[
I(Y;\Z_2)^{\otimes n-2}\otimes I(Z;\Z_2)^{\otimes n-2}\cong \Z_2(y-1)^{\otimes n-2}\otimes \Z_2(z-1)^{\otimes n-2}\cong \Z_2.
\]

Since $\mathbf{c}_n=\mathbf{a}_n+\mathbf{b}_n$, the expression $\Delta_{4,2n-4}(\mathbf{c}_n\times \mathbf{c}_n)$ decomposes as
\[
\begin{aligned}
&\sum\limits_{i=0}^{4}(\mathbf{a}_i\times \mathbf{a}_{4-i})\otimes (\mathbf{a}_{n-i}\times \mathbf{a}_{n-4+i})
+\sum\limits_{i=0}^{4}(\mathbf{a}_i\times \mathbf{b}_{4-i})\otimes (\mathbf{a}_{n-i}\times \mathbf{b}_{n-4+i})\\
&\qquad +\sum\limits_{i=0}^{4}(\mathbf{b}_i\times \mathbf{a}_{4-i})\otimes (\mathbf{b}_{n-i}\times \mathbf{a}_{n-4+i})
+\sum\limits_{i=0}^{4}(\mathbf{b}_i\times \mathbf{b}_{4-i})\otimes (\mathbf{b}_{n-i}\times \mathbf{b}_{n-4+i})
\end{aligned}
\]
Now, we can check that, among the 
right-hand components, the only terms on which the projection of ${\mathfrak v}^{2n-4}$ on $I(Y;\Z_2)^{\otimes n-2}\otimes I(Z;\Z_2)^{\otimes n-2}$ does not vanish are
$ \mathbf{a}_{n-2}\times \mathbf{b}_{n-2}$ and $ \mathbf{b}_{n-2}\times \mathbf{a}_{n-2}$ (represented respectively  by $EZ(\alpha'_{n-2}\otimes\beta'_{n-2})$ and $EZ(\beta'_{n-2}\otimes \alpha'_{n-2})$). Moreover, in $I(Y;\Z_2)^{\otimes n-2}\otimes I(Z;\Z_2)^{\otimes n-2}$, we have
\[
\mathfrak{v}^{2n-4}(\mathbf{a}_{n-2}\times \mathbf{b}_{n-2})=\mathfrak{v}^{2n-4}(\mathbf{b}_{n-2}\times \mathbf{a}_{n-2})=(y-1)^{\otimes n-2}\otimes (z-1)^{\otimes n-2}.
\]

\noindent Consequently, in $I(D;\Z_2)^{\otimes 4}_{D\times D}\otimes  I(Y;\Z_2)^{\otimes n-2}\otimes I(Z;\Z_2)^{\otimes n-2}$, we have
\[\mathfrak{v}^{2n}(\mathbf{c}_n\times \mathbf{c}_n)=\mathfrak{v}^4(\mathbf{a}_2\times \mathbf{b}_2+\mathbf{b}_2\times \mathbf{a}_2)\otimes (y-1)^{\otimes n-2}\otimes (z-1)^{\otimes n-2}.
\]
It remains to check that $\mathfrak{v}^4(\mathbf{a}_2\times \mathbf{b}_2+\mathbf{b}_2\times \mathbf{a}_2)$ does not vanish in $I(D;\Z_2)^{\otimes 4}_{D\times D}$. In order to see this, recall the expression (in $I(D;\Z_2)^{\otimes 4}$) of $\nu^ 4(EZ(\alpha'_2 \otimes \beta'_2))$ which has been obtained in (\ref{expressionexample}). 

By considering, as in \cite{CV}, the projection 
$$I(D;\Z_2)^ {\otimes 4} \to I(D;\Z_2)\otimes \wedge^{ 3}(I(D;\Z_2)) $$
together with the relation $xyx=\bar y$, we see that the expression (\ref{expressionexample}) reduces to
\[(y{x}-x)\otimes (1-y{x}) \wedge (1-\bar{y}) \wedge (1-{x}). \]
Similarly, calculating the image of $\nu^4(EZ(\beta'_2\otimes \alpha'_2))$ in $I(D;\Z_2)\otimes \wedge^{ 3}(I(D;\Z_2))$ yields
\[(y{x}-x)\otimes (1-y{x}) \wedge (1-{y}) \wedge (1-{x}). \]
As in \cite{CV}, we then send the first component to $I(Y;\Z_2)\cong \Z_2$ (through $x\mapsto 1$) and the statement follows from the fact that the sum of the two elements above is the element 
\[
s=(x-1)\wedge (yx-1)\wedge (y-\bar{y})\in \hbox{\text{$\bigwedge^3$}} I(D;\Z_2).
\]
which is shown to be nonzero in \cite[\S3.3.2]{CV}.
\hfill $\Box$

\section*{\bf Acknowledgment} The first author thanks S.~Carter, A.~Dranishnikov, and S.~Ferry for inviting him to participate in the Geometric Topology session at the 52nd Spring Topology and Dynamical Systems Conference.

\newcommand{\arxiv}[1]{{\texttt{\href{http://arxiv.org/abs/#1}{{arXiv:#1}}}}}

\newcommand{\MRh}[1]{\href{http://www.ams.org/mathscinet-getitem?mr=#1}{MR#1}}

\bibliographystyle{plain}

\end{document}